\documentclass[12pt]{article}
\usepackage{graphicx}
\usepackage{amssymb}
\usepackage{amsmath}
\numberwithin{equation}{section}

\begin{document}
\author{Andrey Tydnyuk}
\date{27 March 2007}
\textbf{Explicit Rational Solution of the KZ Equation (example)}

\begin{center} Andrey Tydnyuk \end{center}
E-mail address:andrey.tydniouk@verizon.net\\ 735 Crawford
Ave.Brooklyn, NY 11223. USA. \\

\begin{center}{Abstract}\end{center}

We investigate the Knizhnik-Zamolodchikov linear differential
system. The coefficients of this system are rational functions. We
have proved that the solution of the KZ system is rational when k
is equal to two and n is equal to three (see [5]) . In this paper,
we construct the corresponding solution in the explicit form.
\newpage

\begin{center}{Introduction}\end{center}

We will consider the system of the form:
\begin{equation}\frac{dW}{dz}=-2A(z)W,\end{equation}

where $A(z)$ and $W(z)$ are $3{\times}3$ matrices,
$z_{1}{\ne}z_{2}$. We suppose that $A(z)$ has the form
\begin{equation}A(z)=\frac{P_{1}}{z-z_{1}}+\frac{P_{2}}{z-z_{2}}.\end{equation}

Here:
\begin{equation}P_{1}=\left[\begin{array}{ccc}
  0 & 1 & 0 \\
  1 & 0 & 0 \\
  0 & 0 & 1 \\
\end{array}\right]\end{equation}

\begin{equation}P_{2}=\left[\begin{array}{ccc}
  0 & 0 & 1 \\
  0 & 1 & 0 \\
  1 & 0 & 0 \\
\end{array}\right]\end{equation}

The matrices $P_{1}$ and $P_{2}$ are connected with the matrix
representation of the symmetric group. System (0.1) is a special
case of the Knizhnik-Zamolodchikov [1] , [2]. We have proved that
the solution of system (0.1) is rational [5]. In this paper, we
construct the corresponding solution in the explicit form. We
consider the case when $S_{3}$ and use the method of L. Sakhnovich
[3].

\section{Main Notions, The Coefficients of the solution in the neighborhood of $z = \infty$}

The solution $W(z)$ of system (0.1) has the form [5]:
\begin{equation}
W(z) =
\frac{L_1}{(z-z_1)^2}+\frac{L_2}{(z-z_1)}+\frac{L_3}{(z-z_2)^2}+
\frac{L_4}{(z-z_2)}+z^2G_{-2}+zG_{-1}+G_0.
\end{equation}
In a neighborhood of $z = \infty$ the solution $W(z)$ can be
represented in the form
\begin{equation}
W(z) = \sum_{k = -2}^{\infty}z^{-k}G_k,
\end{equation}
where the coefficients $G_k$ are defined by the relations ( see
[3]).
\begin{equation}
[(q+1)I_3-2T]G_{q+1} = 2\sum_{r + s = q}T_rG_s,\quad r {\geq} 0
\end{equation}
and
\begin{equation}
T_{r} = z_{1}^{r+1}P_{1}+z_{2}^{r+1}P_{2},\quad T = P_{1}+P_{2}.
\end{equation}
The eigenvalues of T are
\begin{equation}
\lambda_{1} = 2, \quad \lambda_{2} = 1, \quad \lambda_{3} = -1.
\end{equation}
The corresponding eigenvectors have the forms:
\begin{equation}
\ell_{1} = \left[
             \begin{array}{c}
               1 \\
               1 \\
               1 \\
             \end{array}
           \right],
\quad \ell_{2} = \left[
             \begin{array}{c}
                0 \\
                1 \\
               -1 \\
             \end{array}
           \right],
\quad \ell_{3} = \left[
             \begin{array}{c}
                2 \\
               -1 \\
               -1 \\
             \end{array}
           \right].
\end{equation}
First, we will begin with finding all of the coefficients from
$G_{-2}$ to $G_{-4}$.\\ The eigenvalues of matrix $2T$ are twice
the eigenvalues of the matrix $T$. Thus we get:
\begin{equation}
\mu_{1} = 4, \quad \mu_{2} = 2, \quad \mu_{3} = -2.
\end{equation}
The eigenvectors remain unchanged.\\ The smallest eigenvalue of
$2T$ is equal to $(-2)$. That is why we begin with the coefficient
$G_{-2}$. From equation (1.3) we can say that
\begin{equation}
(-2I_{3}-2T)G_{-2}=0.
\end{equation}
Using equation (1.6) and (1.8) we conclude that
\begin{equation}
G_{-2} = \left[
             \begin{array}{c}
               2 \\
               -1 \\
               -1 \\
             \end{array}
           \right]
\end{equation}
When the coefficient is $G_{-1}$, equation (1.3) takes the form
\begin{equation}
(-I_{3}-2T)G_{-1}=2T_{0}G_{-2}
\end{equation}
From the last relation we find that \\
\begin{equation}
G_{-1} = 2\left[
             \begin{array}{c}
               -(z_{1}+z_{2}) \\
               z_{2} \\
               z_{1} \\
             \end{array}
           \right].
\end{equation}
When $q + 1 = 0$, we get the relation:
\begin{equation}
-2TG_{0}=2(T_{0}G_{-1}+T_{1}G_{-2}).
\end{equation}
From this we find that\\
\begin{equation}
G_{0} = \left[
             \begin{array}{c}
               -z_{1}^2 + 4z_{1}z_{2}-z_{2}^2 \\
               z_{1}(z_{1}-2z_{2}) \\
               z_{2}(-2z_{1}+z_{2}) \\
             \end{array}
           \right].
\end{equation}
When $q + 1 = 1$ we obtain:
\begin{equation}
(I_{3}-2T)G_{1}=2(T_{0}G_{0}+T_{1}G_{-1}+T_{2}G_{-2}).
\end{equation}
Now we have\\
\begin{equation}
G_{1} = 2\left[
             \begin{array}{c}
               0 \\
               (z_{1}-z_{2})^3 \\
               -(z_{1}-z_{2})^3 \\
             \end{array}
           \right].
\end{equation}
When $q+1=2$:
\begin{equation}
(2I_{3}-2T)G_{2}=2(T_{0}G_{1}+T_{1}G_{0}+T_{2}G_{-1}+T_{3}G_{-2}).
\end{equation}
\textbf{Remark 1.1}\\ When $q+1=-1,0,1$ the matrices
$(q+1)I_{3}-2T$ are invertible. That is why $G_{-1}$,$G_{0}$,and
$G_{1}$ are correctly defined by formulas (1.11),(1.13), and
(1.15).\\ The situation changes when $q+1=2$ because $2$ is an
eigenvalue of the matrix $2T$. In this case, the matrix
$2I_{3}-2T$ is not invertible.\\ The right-hand side of equation
(1.16) has the form:\\
\begin{equation}
            \left[
             \begin{array}{c}
               4(z_{1}-z_{2})^4 \\
               -2(z_{1}-z_{2})^4 \\
               -2(z_{1}-z_{2})^4 \\
             \end{array}
           \right].
\end{equation}
The eigenvalues of $(2I_{3}-2T)$ are
\begin{equation}
\mu_{1}=-2, \quad \mu_{2} = 0, \quad \mu_{3} = 4.
\end{equation}
The right side of (1.16) is the linear combination of the vectors
$\ell_{1}$ and $\ell_{3}$. From relations (1.6), (1.16), and
(1.18) we obtain
\begin{equation}
           G_{2} =  \left[
             \begin{array}{c}
               (z_{1}-z_{2})^4 \\
               -\frac{1} {2}(z_{1}-z_{2})^4 \\
               -\frac{1} {2}(z_{1}-z_{2})^4 \\
             \end{array}
           \right].
\end{equation}
When $q + 1 = 3$ we obtain:
\begin{equation}
(3I_{3} -
2T)G_{3}=2(T_{0}G_{2}++T_{1}G_{1}+T_{2}G_{0}+T_{3}G_{-1}+T_{4}G_{-2}).
\end{equation}

Using our previous results we find that
\begin{equation}
           G_{3} =  \left[
             \begin{array}{c}
               \frac{3}{5}(z_{1}-z_{2})^4(z_{1}+z_{2}) \\
               \frac{1}{5}(z_{1}-z_{2})^3 (6z_{1}^2-25z_{1}z_{2}+9z_{2}^2) \\
               -\frac{1}{5}(z_{1}-z_{2})^3 (9z_{1}^2-25z_{1}z_{2}+6z_{2}^2)\\
             \end{array}
           \right].
\end{equation}
When $q + 1 = 4$ we use the formula
\begin{equation}
(4I_{3} -
2T)G_{4}=2(T_{0}G_{3}++T_{1}G_{2}+T_{2}G_{1}+T_{3}G_{0}+T_{4}G_{-1}+T_{5}G_{-2}).
\end{equation}
The right side of (1.22) has the form
\begin{equation}
            \left[
             \begin{array}{c}
               \frac{9}{5}(z_{1}-z_{2})^4(3z_{1}^2-4z_{1}z_{2}+3z_{2}^2) \\
               \frac{1}{5}(z_{1}-z_{2})^3(6z_{1}^3-8z_{1}^{2}z_{2}-71z_{1}z_{2}^2+33z_{2}^3) \\
               -\frac{1}{5}(z_{1}-z_{2})^3(33z_{1}^3-71z_{1}^{2}z_{2}-8z_{1}z_{2}^2+6z_{2}^3)\\
             \end{array}
           \right].
\end{equation}
The case when $q + 1 = 4$ is similar to the case when $q + 1 = 2$
(see Remark 1.1). The eigenvalues of $(4I_{3}-2T)$ are
\begin{equation}
\mu_{1}=0, \quad \mu_{2} = 2, \quad \mu_{3} = 6.
\end{equation}
The right side of (1.22) is the linear combination of the vectors
$\ell_{2}$ and $\ell_{3}$. From relations (1.22), (1.23), and
(1.24) we obtain

\begin{equation}
           G_{4} =  \left[
             \begin{array}{c}
               \frac{3}{10}(z_{1}-z_{2})^4(3z_{1}^2-4z_{1}z_{2}+3z_{2}^2) \\
               \frac{1}{10}(z_{1}-z_{2})^3(15z_{1}^3-29z_{1}^{2}z_{2}-50z_{1}z_{2}^2+24z_{2}^3) \\
               -\frac{1}{10}(z_{1}-z_{2})^3(24z_{1}^3-50z_{1}^{2}z_{2}-29z_{1}z_{2}^2+15z_{2}^3)\\
             \end{array}
           \right].
\end{equation}
From (1.1) we wind up with the following system:

\begin{align}
L_{2} + L_{4} = G_{1}\\ L_{1} + L_{2}z_{1} + L_{3} + L_{4}z_{2} =
G_{2}\\ 2L_{1}z_{1} + L_{2}z_{1}^2 + 2L_{3}z_{2} + L_{4}z_{2}^2 =
G_{3} \\ 3L_{1}z_{1}^2 + L_{2}z_{1}^3 + 3L_{3}z_{2}^2 +
L_{4}z_{2}^3 = G_{4}
\end{align}

System (1.26)-(1.29) can be written in the matrix form:
\begin{equation}
SX = Y,
\end{equation}
where
\begin{equation}
S = \left[\begin{array}{cccc}
  0 & I_{3} & 0 & I_{3} \\
  I_{3} & z_{1} & I_{3} & z_{2} \\
  2z_{1}I_{3} & z_{1}^2I_{3} & 2z_{2}I_{3} & z_{2}^2I_{3} \\
  3z_{1}^2I_{3} & z_{1}^3I_{3} & 3z_{2}^2I_{3} & z_{2}^3I_{3} \\
\end{array}\right],
\end{equation}

\begin{equation}
X= \mathrm {col} [ L_{1} , L_{2} , L_{3} , L_{4}],
\end{equation}
\begin{equation}
Y = \mathrm {col} [ G_{1} , G_{2} , G_{3} , G_{4}].
\end{equation}

In equation (1.30) the matrices $S$ and $Y$ are known, but the
matrix $X$ is unknown.

From relation (1.30) we get that \begin{equation}X =
S^{-1}Y\end{equation} where

\begin{equation}
S^{-1} = \left[\begin{array}{cccc}
  -\frac {z_{1}z_{2}^2} {(z_{1} - z_{2}) ^ 2}I_{3} & \frac {z_{2}(2z_{1}+z_{2})} {(z_{1} - z_{2}) ^ 2}I_{3} & -\frac {z_{1}+2z_{2}} {(z_{1} - z_{2}) ^ 2}I_{3} & \frac {1} {(z_{1} - z_{2}) ^ 2}I_{3} \\
  \frac {(3z_{1}-z_{2})z_{2}^2} {(z_{1} - z_{2}) ^ 3}I_{3} & -\frac {6z_{1}z_{2}} {(z_{1} - z_{2}) ^ 3}I_{3} & \frac {3z_{1}+z_{2}} {(z_{1} - z_{2}) ^ 3}I_{3} & \frac {2} {(-z_{1} + z_{2}) ^ 3}I_{3} \\
  -\frac {z_{1}^{2}z_{2}} {(z_{1} - z_{2}) ^ 2}I_{3} & \frac {z_{1}(z_{1}+2z_{2})} {(z_{1} - z_{2}) ^ 2}I_{3} & -\frac {2z_{1}+z_{2}} {(z_{1} - z_{2}) ^ 2}I_{3} & \frac {1} {(z_{1} - z_{2}) ^ 2}I_{3}\\
  \frac {z_{1}^2(z_{1}-3z_{2})} {(z_{1} - z_{2}) ^ 3}I_{3} & \frac {6z_{1}z_{2}} {(z_{1} - z_{2}) ^ 3}I_{3} & -\frac {3z_{1}+z_{2}} {(z_{1} - z_{2}) ^ 3}I_{3} & \frac {2} {(z_{1} - z_{2}) ^ 3}I_{3} \\
\end{array}\right].
\end{equation}

Thus, we find that

\begin{equation}
           L_{1} =  \left[
             \begin{array}{c}
               \frac{1}{10}(3z_{1}-7z_{2})(z_{1}-z_{2})^3 \\
               \frac{1}{10}(3z_{1}-7z_{2})(z_{1}-z_{2})^3 \\
               -\frac{1}{5}(3z_{1}-7z_{2})(z_{1}-z_{2})^3 \\
             \end{array}
           \right],
\end{equation}
\begin{equation}
           L_{2} =  \left[
             \begin{array}{c}
               0 \\
               \frac{1}{5}(3z_{1}-7z_{2})(z_{1}-z_{2})^2 \\
               -\frac{1}{5}(3z_{1}-7z_{2})(z_{1}-z_{2})^2 \\
             \end{array}
           \right],
\end{equation}
\begin{equation}
           L_{3} =  \left[
             \begin{array}{c}
               \frac{1}{10}(7z_{1}-3z_{2})(z_{1}-z_{2})^3 \\
               -\frac{1}{5}(7z_{1}-3z_{2})(z_{1}-z_{2})^3 \\
               \frac{1}{10}(7z_{1}-3z_{2})(z_{1}-z_{2})^3 \\
             \end{array}
           \right],
\end{equation}
and
\begin{equation}
           L_{4} =  \left[
             \begin{array}{c}
               0 \\
               \frac{1}{5}(7z_{1}-3z_{2})(z_{1}-z_{2})^3 \\
               -\frac{1}{5}(7z_{1}-3z_{2})(z_{1}-z_{2})^3 \\
             \end{array}
           \right].
\end{equation}

This way we have proved the following statement:\\
\textbf{Proposition 1} \emph{System (0.1) has the following
solution:}
\begin{equation}
W_{1}(z) =
\frac{L_1}{(z-z_1)^2}+\frac{L_2}{(z-z_1)}+\frac{L_3}{(z-z_2)^2}+
\frac{L_4}{(z-z_2)}+z^2G_{-2}+zG_{-1}+G_0.
\end{equation}

\emph{The matrices $G_{k}$ and $L_{k}$ are defined by the
relations (1.9), (1.1), (1.13), and (1.36) - (1.39).}\\ To find
the next solution to the system (0.1) we consider the case

\begin{equation}
g_{k} = 0 \quad when \quad k < 2
\end{equation}
In this case, we have
\begin{equation}
g_{2} = \mathrm {col} [ 0 , 1, -1].
\end{equation}
From relation (1.3) we get:
\begin{equation}
(3I_{3} - 2T)g_{3} = T_{0}g_{2}.
\end{equation}
From this we find that
\begin{equation}
g_{3} = \frac {1} {5} \left[
             \begin{array}{c}
               z_{1} - z_{2} \\
               2z_{1} + 3z_{2} \\
               -3z_{1} - 2z_{2} \\
             \end{array}
           \right].
\end{equation}

In order to find $g_{4}$ we will use equation (1.3) again;

\begin{equation}
(4I_{3} - 2T)g_{4} = T_{0}g_{3} + T_{1}g_{2}.
\end{equation}

The right side of the equation (1.45) has the form:

\begin{equation}
\frac {1} {5} \left[
             \begin{array}{c}
               7(z_{1} - z_{2})(z_{1} + z_{2}) \\
               z_{1} ^ 2 +z_{1}z_{2} + 8z_{2} ^ 2 \\
               -8z_{1} ^ 2 - z_{1}z_{2} - z_{2} ^ 2 \\
             \end{array}
           \right] = \left[ \begin{array}{c}
               \phi \\
               - \frac {\phi} {2} \\
               - \frac {\phi} {2} \\
             \end{array}
           \right] + \left[ \begin{array}{c}
               0 \\
               \psi + \frac {\phi} {2} \\
               -\psi - \frac {\phi} {2} \\
             \end{array}
           \right],
\end{equation}
where
\begin{equation}
\phi = \frac {7} {5} (z_{1} - z_{2})(z_{1} + z_{2}) \quad , \quad
\psi = \frac {1} {5} z_{1} ^ 2 +z_{1}z_{2} + 8z_{2} ^ 2
\end{equation}

Analogously (1.34) we can write:
\begin{equation}
X = S ^ {-1} Y,
\end{equation}
where

\begin{equation}
X= \mathrm {col} [ M_{1} , M_{2} , M_{3} , M_{4}],
\end{equation}

\begin{equation}
Y = \mathrm {col} [ 0 , g_{2} , g_{3} , g_{4}].
\end{equation}

Thus, we can say that

\begin{equation}
           M_{1} =  \left[
             \begin{array}{c}
               \frac {5z_{1} + 3z_{2}} {10(z_{1}-z_{2})} \\
               \frac {z_{1}z_{2} + 9(-2z_{1} ^ 2 + 2z_{1}z_{2} + z_{2} ^ 2)} {30(z_{1} - z_{2}) ^ 2} \\
               -\frac {z_{1}z_{2} - 3z_{1}(z_{1}-4z_{2})} {30(z_{1} - z_{2}) ^ 2}\\
             \end{array}
           \right],
\end{equation}
\begin{equation}
           M_{2} =  \left[
             \begin{array}{c}
               -\frac {4(z_{1} + z_{2})} {5(z_{1}-z_{2})^2} \\
               -\frac {-24z_{1}^2 + 46z_{1}z_{2} - 12z_{2}^2} {15(z_{1} - z_{2}) ^ 3} \\
               \frac {-12z_{1}^2 + + 46z_{1}z_{2} - 24z_{2}^2} {15(z_{1} - z_{2}) ^ 3} \\
             \end{array}
           \right],
\end{equation}
\begin{equation}
           M_{3} =  \left[
             \begin{array}{c}
               -\frac{3z_{1} + 5z_{2}} {10(z_{1} - z_{2})}\\
               \frac{z_{1}z_{2} + 3(4z_{1} - z_{2})z_{2}}{30(z_{1} - z_{2}) ^ 2} \\
               -\frac{z_{1}z_{2} + 9(z_{1} ^ 2 + 2z_{1}z_{2} - 2z_{2} ^2}{30(z_{1} - z_{2}) ^ 2}\\
             \end{array}
           \right],
\end{equation}
and
\begin{equation}
           M_{4} =  \left[
             \begin{array}{c}
               \frac {4(z_{1} + z_{2})} {5(z_{1}-z_{2})^2} \\
               \frac {-24z_{1}^2 + 46z_{1}z_{2} - 12z_{2}^2} {15(z_{1} - z_{2}) ^ 3} \\
               -\frac {-12z_{1}^2 + 46z_{1}z_{2} - 24z_{2}^2} {15(z_{1} - z_{2}) ^ 3} \\
             \end{array}
           \right].
\end{equation}

This way we have proved the following statement:\\
\textbf{Proposition 2} \emph{System (0.1) has the following
solution:}
\begin{equation}
W_{2}(z) =
\frac{M_1}{(z-z_1)^2}+\frac{M_2}{(z-z_1)}+\frac{M_3}{(z-z_2)^2}+
\frac{M_4}{(z-z_2)}.
\end{equation}

\emph{The matrices $M_{k}$ are defined by the relations (1.51) -
(1.54).}\\

The next solution of system (0.1) has the form
\begin{equation}
W_{3}(z) =
\frac{N_1}{(z-z_1)^2}+\frac{N_2}{(z-z_1)}+\frac{N_3}{(z-z_2)^2}+
\frac{N_4}{(z-z_2)}.
\end{equation}

In order to find $N_{1}$ , $N_{2}$ , $N_{3}$ , and $N_{4}$ we
consider the case when
\begin{equation}
G_{k} = 0 \quad when \quad k < 4 \quad and \quad G_{4} = \ell_{1}.
\end{equation}

From relations (1.34) and (1.35) we have
\begin{equation}
N_{1} = N_{3} = \frac {1} {(z_{1} - Z_{2}) ^ 2} \ell_{1}
\end{equation}
\begin{equation}
N_{2} = -N_{4} = \frac {2} {(-z_{1} + z_{2}) ^ 3}.
\end{equation}

The main theorem follows from Propositions 1 - 3.

\textbf{Theorem 1}\\ \emph{The general solution of system (0.1)
has the form:}
\begin{equation}
W_{z} = \alpha_{1}W_{1}(z) + \alpha_{2}W_{2}(z) +
\alpha_{3}W_{3}(z),
\end{equation}
\emph{where $\alpha_{1}$ , $\alpha_{2}$ , and $\alpha_{3}$ are
arbitrary constants.}
\\
\\
\\

\begin{center}{References}\end{center}
1. Chervov A., Talalaev D., KZ equation, G-opers, quantum
Drinfeld-Sokolov reduction and quantum Cayley-Hamilton identity,\\
arXiv:hep-th/0607250, 2006.\\ 2. Etingof P.I., Frenkel I.B.,
Kirillov A.A. jr.,Lectures on Representation Theory and
Knizhnik-Zamolodchikov Equations, Amer. Math. Society, 1998.\\ 3.
Sakhnovich L.A., Meromorphic Solutions of Linear Differential
Systems, Painleve Type Functions, Journal of Operator and Matrix
Theory, v. 1, 87-111, 2007.\\ 4. Sakhnovich L.A., Rational
Solution of the KZ Equation, existence and construction,
arXiv:math.ph/0609067, 2006.\\ 5. Tydnyuk A.I., Rational Solution
of the KZ Equation (example), \\arXiv:math.CA/0612153, 2006.

\end{document}